\title{Automorphism groups of Beauville surfaces}

\author{Gareth A. Jones\\
School of Mathematics\\
University of Southampton\\
Southampton SO17  1BJ, U.K.\\
{\tt G.A.Jones@maths.soton.ac.uk}
}

\documentclass[12pt]{article}
\usepackage[latin1]{inputenc}
\usepackage{a4wide}
\usepackage{latexsym}
\usepackage{amsmath}
\usepackage{amssymb}
\usepackage{amsfonts}
\newtheorem{thm}{Theorem}[section]
\newtheorem{lemma}[thm]{Lemma}
\newtheorem{cor}[thm]{Corollary}
\newtheorem{prop}[thm]{Proposition}
\date{}
\begin{document} 

\maketitle

\begin{abstract}
A Beauville surface of unmixed type is a complex algebraic surface which is the quotient of the product of two curves of genus at least 2 by a finite group $G$ acting freely on the product, where $G$ preserves the two curves and their quotients by $G$ are isomorphic to the projective line, ramified over three points. We show that the automorphism group $A$ of such a surface has an abelian normal subgroup $I$ isomorphic to the centre of $G$, induced by pairs of elements of $G$ acting compatibly on the curves (a result obtained independently by Fuertes and Gonz\'alez-Diez). Results of Singerman on inclusions between triangle groups imply that $A/I$ is isomorphic to a subgroup of the wreath product $S_3\wr S_2$, so $A$ is a finite solvable group. Using constructions based on Lucchini's work on generators of special linear groups, we show that every finite abelian group can arise as $I$, even if one restricts the index $|A:I|$ to the extreme values 1 or 72.
\end{abstract}


\section{Introduction}

A Beauville surface $\cal S$ is a complex algebraic surface which is rigid and is isogenous to a higher product, that is, it has the form $({\cal C}_1\times{\cal C}_2)/G$, where ${\cal C}_1$ and ${\cal C}_2$ are complex algebraic curves of genus at least $2$, and $G$ is a finite group acting freely on their product. Here we restrict attention to Beauville surfaces of unmixed type, where $G$ preserves the factors ${\cal C}_i$, and rigidity means that ${\cal C}_i/G$ is isomorphic to the projective line ${\mathbb P}^1({\mathbb C})$, with the covering ${\cal C}_i\to{\mathbb P}^1({\mathbb C})$ ramified over three points. (This implies that each ${\cal C}_i$ carries a regular dessin, in the sense of Grothendieck's theory of {\it dessins d'enfants\/}~\cite{Gro}, so by Bely\u\i's Theorem~\cite{Bel} $\cal S$ is defined over an algebraic number field.) The first example, with ${\cal C}_1$ and ${\cal C}_2$ Fermat curves, was introduced by Beauville in~\cite{Bea}, and subsequently these surfaces have been studied by geometers such as Bauer, Catanese and Grunewald~\cite{BCG05, BCG06, Cat}. Much recent research has concentrated on the question of which groups $G$, and in particular which simple groups, can arise in this context~\cite{Fai, FMP, FG, FJ, GP, GLL, GM}. Here we take a different point of view, studying the automorphism group $A$ of $\cal S$, and how it is related to $G$. 

The automorphisms of a Beauville surface $\cal S$ are of two types, which we will call {\sl direct\/} or {\sl indirect\/} as they lift to automorphisms of ${\cal C}_1\times{\cal C}_2$ which either preserve or transpose the two factors ${\cal C}_i$. The direct automorphisms form a subgroup $A^0={\rm Aut}^0{\cal S}$ of index at most $2$ in the group $A={\rm Aut}\,{\cal S}$ of all automorphisms of $\cal S$. Indirect automorphisms exist if and only if the two curves ${\cal C}_i$ are isomorphic, with the actions of $G$ on them transposed by an automorphism of $G$. We will concentrate mainly on the group $A^0$. We show that this has an abelian normal subgroup $I={\rm Inn}\,{\cal S}$ isomorphic to the centre $Z(G)$ of $G$; its elements, which we call the {\sl inner automorphisms\/} of $\cal S$, are induced by pairs of elements of $G$ acting compatibly on the curves ${\cal C}_i$. In most cases $A^0=I$, but in some cases either curve ${\cal C}_i$ may have additional automorphisms corresponding to non-identity permutations of the three ramification points, and these may induce what we call {\sl outer automorphisms\/} of $\cal S$. The ramification condition on each ${\cal C}_i$ is equivalent to $G$ arising as a quotient of a hyperbolic triangle group by a normal surface subgroup uniformising ${\cal C}_i$. Using results of Singerman~\cite{Sin2} on inclusions between triangle groups, we show that the {\sl outer automorphism group\/} ${\rm Out}\,{\cal S}=A/I$ is isomorphic to a subgroup of the wreath product $S_3\wr S_2$, a semidirect product of $S_3\times S_3$ by $S_2$, so that $A$ is a finite solvable group. A result of Lucchini~\cite{Luc} on generators of special linear groups allows us to show that every finite abelian group is isomorphic to ${\rm Inn}\,{\cal S}$ for some Beauville surface $\cal S$; moreover, we show that $\cal S$ can be chosen here so that ${\rm Out}\,{\cal S}$ is as large or as small as possible, namely isomorphic to $S_3\wr S_2$ or the trivial group. We also give further examples with ${\rm Out}\,{\cal S}$ lying between these two extremes, for instance showing that every finite generalised dihedral group can arise as ${\rm Aut}\,{\cal S}$ for some $\cal S$.  

The author thanks Fabrizio Catanese, Yolanda Fuertes, Gabino Gonz\'alez-Diez, David Torres-Teigell and J\"urgen Wolfart for helpful comments on earlier drafts of this paper.

\section{Beauville surfaces and structures}

A {\sl Beauville surface\/} (of unmixed type) is a compact complex surface $\cal S$ such that
\begin{description}
\item[\rm(a)] $\cal S$ is isogenous to a higher product, that is, ${\cal S}\cong ({\cal C}_1\times {\cal C}_2)/G$ where ${\cal C}_1$ and ${\cal C}_2$ are algebraic curves of genus at least $2$ and $G$ is a finite group acting faithfully on ${\cal C}_1$ and ${\cal C}_2$ by holomorphic transformations in such a way that it acts freely on ${\cal C}_1\times{\cal C}_2$;
\item[\rm(b)] each ${\cal C}_i/G$ is isomorphic to the projective line ${\mathbb P}^1({\mathbb C})$, and the covering ${\cal C}_i\to {\cal C}_i/G$ is ramified over three points.
\end{description}
(We will not consider Beauville surfaces of mixed type, where $G$ contains elements which transpose the two curves ${\cal C}_i$.) Condition~(b) is equivalent to each curve ${\cal C}_i$ admitting a regular {\it dessin\/} in the sense of Grothendieck's theory of {\it dessins d'enfants\/}~\cite{Gro}, or equivalently an orientably regular hypermap~\cite{JS}, with $G$ acting as the orientation-preserving automorphism group. By Bely\u\i's Theorem~\cite{Bel}, this implies that the curves ${\cal C}_i$ and the coverings ${\cal C}_i\to{\mathbb P}^1({\mathbb C})$ are defined over algebraic number fields, and hence the same applies to $\cal S$.

A finite group $G$ arises in this way if and only if it has generating triples $a_i, b_i$ and $c_i$ for $i=1, 2$, of orders $l_i, m_i$ and $n_i$, such that
\begin{description}
\item[\rm(1)] $a_ib_ic_i=1$ for each $i=1, 2$,
\item[\rm(2)] $l_i^{-1}+m_i^{-1}+n_i^{-1}<1$ for each $i=1, 2$, and
\item[\rm(3)] no non-identity power of $a_1, b_1$ or $c_1$ is conjugate in $G$ to a power of $a_2, b_2$ or $c_2$.
\end{description}

\noindent We will call such a pair of triples $(a_i, b_i, c_i)$ a {\sl Beauville structure\/} for $G$. Property~(1) is equivalent to $G$ being a smooth quotient $\Delta_i/K_i$ of a triangle group $\Delta_i=\Delta(l_i, m_i, n_i)$ by a normal surface subgroup $K_i$ uniformising ${\cal C}_i$, with $a_i$, $b_i$ and $c_i$ the local monodromy permutations for the covering ${\cal C}_i\to {\cal C}_i/G$ at the three ramification points; we call $l_i, m_i$ and $n_i$ the {\sl elliptic periods\/} of $\Delta_i$. Property~(2) is equivalent to each ${\cal C}_i$ having genus at least $2$, so that $\Delta_i$ acts on the hyperbolic plane $\mathbb H$, with ${\cal C}_i\cong{\mathbb H}/K_i$, and property~(3) is equivalent to $G$ acting freely on  ${\cal C}_1\times{\cal C}_2$.

To show that a pair such as $a_1$ and $a_2$ satisfy property~(3) it is sufficient to verify that, for each prime $p$ dividing both $l_1$ and $l_2$, the element $a_1^{l_1/p}$ of order $p$ is not conjugate to any of the elements $a_2^{kl_2/p}$ where $k=1, 2, \ldots, p-1$; in particular, if $l_1$ is prime then it is sufficient to verify this for $a_1$.  Similar remarks apply to every other pair chosen from the two triples, so in particular property~(3) is satisfied if $l_1m_1n_1$ is coprime to $l_2m_2n_2$.

A {\sl triple\/} in a group $G$ will always mean an ordered triple $(a, b, c)$ of elements of $G$ such that $abc=1$; it is a {\sl generating triple\/} if $a, b$ and $c$ (and hence any two of them) generate $G$, and it has {\sl type\/} $(l, m, n)$ if these {\sl periods\/} $l, m$ and $n$ are the orders of $a, b$ and $c$. Two types are {\sl equivalent\/} if they differ by a permutation of their periods, so that they define the same triangle groups, but with different generating triples. We will say that a Beauville structure, in the notation above, has {\sl type\/} $(l_1, m_1, n_1; l_2, m_2, n_2)$; here equivalence allows us to permute or transpose the two multisets $\{l_i, m_i, n_i\}$.

\section{Automorphism groups}

Let ${\cal S}=({\cal C}_1\times{\cal C}_2)/G$ be a Beauville surface, as described in Section~2. Any automorphism $\alpha$ of $\cal S$ lifts to an automorphism $\overline\alpha$ of ${\cal C}_1\times{\cal C}_2$, and this is of one of the two following types, where $A_i:={\rm Aut}\,{\cal C}_i$ for $i=1, 2$ (see~\cite{Cat}):

\begin{enumerate}
\item $\overline\alpha=(\alpha_1, \alpha_2):p=(p_1, p_2)\mapsto (p_1\alpha_1, p_2\alpha_2)$ where $\alpha_i\in A_i$ for $i=1, 2$, or
\item $\overline\alpha=(\phi_1, \phi_2):p=(p_1, p_2)\mapsto (p_2\phi_2, p_1\phi_1)$ where $\phi_1:{\cal C}_1\to{\cal C}_2$ and $\phi_2:{\cal C}_2\to{\cal C}_1$ are isomorphisms.
\end{enumerate}
Let us call automorphisms $\alpha$ and $\overline\alpha$ {\sl direct\/} or {\sl indirect\/}, as $\overline\alpha$ is of type~1 or 2 respectively. The direct automorphisms form a subgroup of ${\rm Aut}\,({\cal C}_1\times{\cal C}_2)$ isomorphic to $A_1\times A_2$. This is the whole of ${\rm Aut}\,({\cal C}_1\times{\cal C}_2)$ unless ${\cal C}_1\cong{\cal C}_2$ ($\cong{\cal C}$, say), in which case ${\rm Aut}\,({\cal C}_1\times{\cal C}_2)$ is isomorphic to the wreath product ${\rm Aut}\,{\cal C}\wr S_2$, a semidirect product of $A_1\times A_2\cong ({\rm Aut}\,{\cal C})^2$ by a complement $S_2$ transposing the direct factors.

\subsection{Direct automorphisms}

A direct automorphism $(\alpha_1,\alpha_2)\in A_1\times A_2$ of ${\cal C}_1\times{\cal C}_2$  induces an automorphism of $\cal S$ if and only if, whenever $p=(p_1,p_2)\in{\cal C}_1\times{\cal C}_2$ and $q=(q_1,q_2)\in{\cal C}_1\times{\cal C}_2$ are equivalent under $G$, then so are their images under $(\alpha_1,\alpha_2)$. More explicitly, we require that if $p_ig=q_i$ for $g\in G$ and $i=1, 2$, then there is some $h\in G$ such that $(p_i\alpha_i)h=q_i\alpha_i$ for $i=1, 2$. If it exists, then such an element $h$ is unique, and independent of the points $p_i$, since a generic point has trivial stabiliser, so we can write this condition as $\alpha_i^{-1}g\alpha_i=h\in G$ for $i=1, 2$. We require this for all $g\in G$, so $\alpha_i$ must be an element of the normaliser $N_i:=N_{A_i}(G)$ of $G$ in $A_i$ for $i=1, 2$, with the induced automorphisms $\beta_i: g\mapsto h$ of $G$ satisfying $\beta_1=\beta_2$, i.e.~the natural homomorphisms $\theta_i:N_i\to {\rm Aut}\,G$ induce
\[\theta=(\theta_1, \theta_2): N_1\times N_2\to{\rm Aut}\,G\times{\rm Aut}\,G\]
sending $(\alpha_1,\alpha_2)$ to an element of the diagonal subgroup $E$ of ${\rm Aut}\,G\times{\rm Aut}\,G$. Thus $(\alpha_1, \alpha_2)$ induces an automorphism of $\cal S$ if and only if it lies in $N:=\theta^{-1}(E)$.

Such a pair $(\alpha_1, \alpha_2)$ acts trivially on $\cal S$ if and only if there is some $g\in G$ such that $p_i\alpha_ig=p_i$ for all $p\in{\cal C}_i$ ($i=1, 2$), i.e.~$\alpha_i\in G$ for $i=1, 2$. Thus the kernel of the action of $N$ on $\cal S$ is the diagonal subgroup $D$ of $G\times G\leq N_1\times N_2\leq A_1\times A_2$, so the group $A^0={\rm Aut}^0{\cal S}$ of direct automorphisms of $\cal S$ has the form
\begin{equation}
A^0\cong N/D.
\end{equation}

In all cases, $G\leq N_i$ and the restriction $\theta_i|_G:G\to{\rm Aut}\,G$ is simply the action of $G$ on itself by conjugation, with kernel equal to the centre $Z:=Z(G)$, so $N$ contains a normal subgroup
\[M=N\cap(G\times G)=\{(g_1, g_2)\in G\times G\mid g_1g_2^{-1}\in Z\}=D\times Z,\]
where the direct factor $Z$ can be taken to be the centre of either direct factor of $G\times G$. Hence $A^0$ contains a normal subgroup $I:={\rm Inn}\,{\cal S}\cong M/D\cong Z$, consisting of the {\sl inner automorphisms\/} of $\cal S$, those induced by compatible pairs of elements of $G$ acting on the curves ${\cal C}_i$, or equivalently by elements of $Z$ acting naturally on one curve and fixing the other. Since $I$ is isomorphic to the centre $Z$ of $G$, it is finite and abelian. We will call the quotient group $A^0/I\cong N/M$ the {\sl direct outer automorphism group\/} ${\rm Out}^0{\cal S}$ of $\cal S$.

In most cases $G=A_i$, so $G=N_i$ and $N=M=D\times Z$; then~(1) implies that in such cases we have
\begin{equation}
A^0=I\cong Z.
\end{equation}
The only possible exceptions to~(2) are where $G$ is a proper subgroup of $N_i$ for some $i$, so that $\Delta_i$ is a proper normal subgroup of a Fuchsian group $\tilde\Delta_i$, with $\tilde\Delta_i/K_i\cong N_i$. As shown by Singerman in~\cite{Sin2}, since $\tilde\Delta_i$ contains a triangle group it must also be a triangle group. He showed that the only possibilities for a proper normal inclusion $\Delta\triangleleft\tilde\Delta$ of one hyperbolic triangle group $\Delta=\Delta_i$ in another triangle group $\tilde\Delta=\tilde\Delta_i$ are (up to permutations of the periods) of the form
\[{\rm (a)} \; \Delta(s, s, t)\triangleleft\Delta(2, s, 2t),\quad
{\rm (b)}\; \Delta(t, t, t)\triangleleft\Delta(3, 3, t),\quad
{\rm (c)}\; \Delta(t, t, t)\triangleleft\Delta(2, 3, 2t)\]
for some integers $s$ and $t$, with $\tilde\Delta/\Delta\cong C_2$, $C_3$ or $S_3$ respectively. In all three cases, at least two of the three periods of $\Delta$ are equal, so we have:

\begin{prop} If a Beauville structure on a group $G$ has type $(l_1, m_1, n_1; l_2, m_2, n_2)$, and for each $i$ the periods $l_i, m_i$ and $n_i$ are mutually distinct, then the direct automorphism group ${\rm Aut}^0{\cal S}$ of the corresponding Beauville surface $\cal S$ is isomorphic to the centre $Z(G)$ of $G$. \hfill$\square$
\end{prop}

In each of the exceptional cases (a), (b) and (c), let $u, v$ and $w$ be the canonical elliptic generators of $\tilde\Delta$, with $uvw=1$. In case~(a) we can take $\Delta$ to have elliptic generators $v^u=uvu$, $v$ and $w^2$, of orders $s, s$ and $t$, with $v^u.v.w^2=(uv)^2w^2=1$; then conjugation by $u$ induces an automorphism of $\Delta$ transposing $v^u$ and $v$. In case (b) we can take $\Delta$ to have generators $w$, $w^{u^2}=uwu^{-1}$ and $w^u=u^{-1}wu$, all of order $t$, with $w.w^{u^2}.w^u=(wu)^3=v^{-3}=1$; conjugation by $u$ induces a $3$-cycle on these generators. The inclusion in case~(c) is a composition of two normal inclusions
\[\Delta(t, t, t)\triangleleft\Delta(3, 3, t)\triangleleft\Delta(2, 3, 2t),\]
of types (b) and (a), with generators $v^u$, $v$ and $w^2$ for $\Delta(3,3,t)$ and hence $w^2$, $(w^2)^{(v^u)^2}$ and $(w^2)^{v^u}$ for $\Delta(t,t,t)$; in this case
\[w^2.(w^2)^{(v^u)^2}.(w^2)^{v^u}=(w^2uvu)^3=(wu)^3=v^{-3}=1,\]
and the element $v^u$ induces a $3$-cycle on these generators. In each of cases~(a), (b) and (c), if a normal subgroup $K$ of $\Delta$ is also normal in $\tilde\Delta$, then the generators $\alpha$ of $N=\tilde\Delta/K$ corresponding to the elliptic generators of $\tilde\Delta$ must induce automorphisms $\beta$ of $G=\Delta/K$, acting as above on the generators of $G$ corresponding to those of $\Delta$. Conversely, if a quotient $G=\Delta/K$ of $\Delta$ has automorphisms $\beta$ acting in this way on its generators, then by forming appropriate semidirect products with the groups $\langle\beta\rangle$ (twice in case~(c)), we see that $K$ is normal in $\tilde\Delta$ and $G$ is a normal subgroup of index $2, 3$ or $6$ in $N=\tilde\Delta/K$.

If $G$ admits a Beauville structure then by applying the above arguments to $G$ as a quotient of the appropriate triangle groups $\Delta_i$, one can determine the groups $N_i$. By considering the automorphisms of $G$ induced by each $N_i$ one can then determine $N$ and hence $A^0$, using~$(1)$.

For each $i=1, 2$ we have $N_i/G\cong\tilde\Delta_i/\Delta_i$, isomorphic to a subgroup of $S_3$, so $A^0/I\cong N/M$ is isomorphic to a subgroup of $S_3\times S_3$. Since this group has order $36$, and is solvable of derived length $2$, we have the following generalisation of Proposition~3.1, most of which has also been obtained by Fuertes and Gonz\'alez-Diez in~\cite{FG2}:

\begin{prop} If $\cal S$ is a Beauville surface, obtained from a Beauville structure on a group $G$, then the direct automorphism group ${\rm Aut}^0{\cal S}$ of $\cal S$ has an abelian normal subgroup $I\cong Z(G)$ with quotient group ${\rm Out}^0{\cal S}$ isomorphic to a subgroup of $S_3\times S_3$. In particular ${\rm Aut}^0{\cal S}$ is solvable, of derived length at most $3$, and it has order dividing $36\,|Z(G)|$. \hfill$\square$
\end{prop}

See Section~4.5 and Theorem~5.6, and also~\cite[\S 5]{FG2}, for instances in which  $|{\rm Out}^0{\cal S}|$ attains the upper bound of $36$.

There is a natural interpretation of the embedding of ${\rm Out}^0{\cal S}$ in $S_3\times S_3$. The covering ${\cal C}_i\to{\cal C}_i/G\cong{\mathbb H}/\Delta_i\cong {\mathbb P}^1({\mathbb C})$ is a Bely\u\i\/ function, that is, a non-constant meromorphic function unbranched outside three points; it is convenient to apply a M\"obius transformation of ${\mathbb P}^1({\mathbb C})$ so that these points are $0, 1$ and $\infty$. They represent the orbits of $G$ on ${\cal C}_i$ with nontrivial stabilisers, namely cyclic groups of orders $l_i, m_i$ and $n_i$, so in the language of orbifolds they are cone-points of these orders. Since $\Delta_i$ is normal in $\tilde\Delta_i$, there is an action of $\tilde\Delta_i$ as a group of automorphisms of ${\mathbb P}^1({\mathbb C})$ leaving invariant the set $B=\{0, 1, \infty\}$. The kernel of this action is $\Delta_i$, and there is a corresponding action of $N_i\cong\tilde\Delta/K_i$, with kernel $G\cong\Delta_i/K$, inducing a subgroup of $S_3$ on $B$. We therefore obtain a product action of $N_1\times N_2$, with kernel $G\times G$, on $B^2$; this induces a subgroup of $S_3\times S_3$, preserving the equivalence relations $\equiv_i$ on $B^2$ defined by pairs having the same $i$th component, for $i=1, 2$. Restricting this action to $N$, with kernel $M$, we obtain a faithful action of $N/M\cong {\rm Out}^0{\cal S}$ on $B^2$ as a subgroup of $S_3\times S_3$.

\subsection{Splitting properties}

In the normal inclusions (a) and (b), the triangle group $\tilde\Delta$ is a semidirect product of $\Delta$ by $C_2$ or $C_3$ respectively: in case (a) we can take $\langle u\rangle$ as a complement for $\Delta$ in $\tilde\Delta$, and in case~(b) we can take $\langle u\rangle$ or $\langle v\rangle$. It follows that if $N_i$ contains $G$ with index $2$ or $3$ then it is a split extension of $G$, with the image of $u$ or $v$ generating the complement. However, in case~(c) $\tilde\Delta$ does not split over $\Delta$: a finite subgroup of a hyperbolic (or euclidean) triangle group must be cyclic, so $\tilde\Delta$ has no subgroups isomorphic to the quotient group $\tilde\Delta/\Delta\cong S_3$.

We have $A^0/I\cong N/M$, and it follows from the above argument that if this group has order $2$ or $3$ then $A^0$ and $N$ split over $I$ and $M$. If $A^0/I$ has order $4$ then involutions in $N_i\setminus G$ ($i=1, 2$) commute and generate a Klein four-group $V_4$ complementing $M$ in $N$, so $A^0$ is a semidirect product of $I$ by $V_4$. A similar argument applies if $|A^0/I|=9$, giving a complement $C_3\times C_3$.

These splitting properties imply that not every group which satisfies the conclusions of Proposition~3.2 can arise as the direct automorphism group of a Beauville surface.

\medskip

\noindent{\bf Example 3.1} For each integer $e\geq 2$ the generalised quaternion group
\[Q=\langle g, h\mid g^{2^e}=1,\; g^h=g^{-1},\; h^2=g^{2^{e-1}}\rangle\]
of order $2^{e+1}$ is an extension of a cyclic normal subgroup $\langle g\rangle$ of order $2^e$ by $C_2$, so it satisfies the conclusions of Proposition~3.2. It is nonabelian, so if it arises as the direct automorphism group $A^0$ of a Beauville surface then $A^0>I$, and hence $|A^0:I|=2$ or $4$ since this index is a power of $2$ dividing $36$. Any normal subgroup of index $2$ or $4$ in $Q$ has an abelian quotient, so it contains the commutator $[g,h]=g^{-2}$; it therefore contains $g^{2^{e-1}}\negthinspace$, which is the only element of order $2$ in $Q$, so $Q$ cannot split over such a normal subgroup. It follows that no Beauville surface can satisfy $A^0\cong Q$.

\subsection{Indirect automorphisms}

Any indirect automorphism of ${\cal C}_1\times{\cal C}_2$ has the form
\begin{equation}
t:{\cal C}_1\times{\cal C}_2\to{\cal C}_1\times{\cal C}_2\,,\;
(p_1, p_2)\mapsto(p_2\phi_2, p_1\phi_1),
\end{equation}
where $\phi_1:{\cal C}_1\to{\cal C}_2$ and  $\phi_2:{\cal C}_2\to{\cal C}_1$ are isomorphisms of curves (in which case there are $|{\rm Aut}\,{\cal C}_1|=|{\rm Aut}\,{\cal C}_2|$ possibilities for each). Equivalently, we can identify ${\cal C}_1$ and ${\cal C}_2$ via $\phi_1$, call the resulting curve $\cal C$, and define
\[t:{\cal C}^2\to{\cal C}^2\,,\;(p_1, p_2)\mapsto(p_2\phi, p_1)\]
where $\phi=\phi_2$ is now an automorphism of $\cal C$.

In order to allow for the possibility (in fact, the necessity) of $G$ acting in different ways on the two factors, it is useful when considering indirect automorphisms to regard $G$ as acting on each ${\cal C}_i$ by means of a faithful representation $\rho_i:G\to G_i\leq A_i$. Thus, when we write ${\cal S}=({\cal C}_1\times{\cal C}_2)/G$, we are really factoring out the action of the image of the diagonal subgroup of $G\times G$ under the product representation $\rho_1\times\rho_2:G\times G\to A_1\times A_2$.

Any indirect automorphism $\tau$ of $\cal S$ must be induced by some $t$ of the form~$(3)$. If they exist, such automorphisms of $\cal S$ are all obtained by composing a specific indirect automorphism $\tau$ (in either order) with an arbitrary direct automorphism of $\cal S$. Now $t$ induces an automorphism of $\cal S$ if and only if, whenever some $g\in G$ sends $(p_1, p_2)$ to $(q_1, q_2)$ in ${\cal C}_1\times{\cal C}_2$, there exists some $h\in G$ sending $(p_1, p_2)t=(p_2\phi_2, p_1\phi_1)$ to $(q_1, q_2)t=(q_2\phi_2, q_1\phi_1)$, that is, $p_2\phi_2(h\rho_1)=p_2(g\rho_2)\phi_2$ and $p_1\phi_1(h\rho_2)=p_1(g\rho_1)\phi_1$. If such an element $h$ exists then it is  unique and is independent of $p_1$ and $p_2$, so we require $\phi_2(h\rho_1)=(g\rho_2)\phi_2$ and $\phi_1(h\rho_2)=(g\rho_1)\phi_1$. The mapping $g\mapsto h$, if it exists, is an automorphism $\zeta$ of $G$; we then require $(g\zeta)\rho_1=\phi_2^{-1}(g\rho_2)\phi_2$ and $(g\zeta)\rho_2=\phi_1^{-1}(g\rho_1)\phi_1$ for all $g\in G$, so that $\phi_2$ induces an equivalence between the representations $\rho_2$ and $\zeta\circ\rho_1$ of $G$, while $\phi_1$ induces an equivalence between $\rho_1$ and $\zeta\circ\rho_2$. The two representations $\rho_1$ and $\rho_2$ cannot be equivalent, for otherwise condition~(3) of a Beauville structure would not hold. It follows that $\zeta$ transposes the two equivalence classes of representations $\rho_i$; in particular, $\zeta$ cannot be an inner automorphism of $G$. Conversely, if $G$ has an automorphism transposing these two classes, then the isomorphisms $\phi_i$ realising these equivalences give an indirect automorphism of $\cal S$. Thus we have proved:

\begin{prop}
A Beauville surface ${\cal S}=({\cal C}_1\times{\cal C}_2)/G$ has an indirect automorphism if and only if ${\cal C}_1\cong{\cal C}_2$ and $G$ has an automorphism $\zeta$ transposing the equivalence classes of its representations on ${\cal C}_1$ and ${\cal C}_2$. \hfill$\square$
\end{prop}

\begin{cor}
If a Beauville surface ${\cal S}=({\cal C}_1\times{\cal C}_2)/G$ has an indirect automorphism, then the corresponding Beauville structure for $G$ must consist of two triples of equivalent types.
\end{cor}

\noindent{\sl Proof.} The automorphism $\zeta$ of $G$ must preserve the orders of the stabilisers of points in the representations of $G$ on the two curves. \hfill$\square$

\medskip

We will initially use Corollary~3.4 to show that various Beauville surfaces do not possess indirect automorphisms. Later, in Section~6, we will consider indirect automorphisms in more detail, and give examples of Beauville surfaces with such automorphisms.

\section{Examples of direct automorphism groups}

This section contains some specific examples of automorphism groups which illustrate the general results proved in Section~3.

\subsection{Examples with trivial automorphism groups}

In~\cite{FG}, Fuertes and Gonz\'alez-Diez have shown that the symmetric groups $G=S_n$ admit Beauville structures for all $n\geq 5$. They give examples of Beauville structures of types $(2, n-2, n-3; 2, n, n-1)$ and $(2, 4(n-6), n-2; 2, n-1, n)$ respectively for even and odd $n>10$. These types satisfy the conditions of Proposition~3.1, and the groups $S_n$ have trivial centres, so the corresponding Beauville surfaces have $A^0=1$. In each case the two triples have inequivalent types, so it follows from Corollary~3.4 that $A=1$ also.

\subsection{Examples with automorphism group $C_2$}

Let $G$ be the simple group $L_2(p)=SL_2(p)/\{\pm I\}$ for some prime $p\equiv 1$ mod~$(24)$. Since $p\equiv 1$ mod~$(4)$ there is some $u$ in the field $ {\mathbb F}_p$ such that $u^2=-1$.
Let
\[
a_1=\pm\Big(\,\begin{matrix}1&u\cr 0&1\end{matrix}\,\Big),
\quad
b_1=\pm\Big(\,\begin{matrix}1&0\cr u&1\end{matrix}\,\Big)
\quad{\rm and}\quad
c_1=(a_1b_1)^{-1}=\pm\Big(\,\begin{matrix}1&-u\cr -u&0\end{matrix}\,\Big),
\]
so that $a_1$ and $b_1$ have order $p$, and $c_1$ has order $3$. The maximal subgroups of $G$ have been classified by Dickson~\cite[Ch.~XII]{Dic}, and of these, only the stabilisers of points in the projective line ${\mathbb P}^1(p)$ have order divisible by $p$; since $a_1$ and $b_1$ have no common fixed points no such subgroup contains them both, so they generate $G$. 
Now let
\[
a_2=\pm\Big(\,\begin{matrix}0&v\cr -v^{-1}&w\end{matrix}\,\Big),
\;
b_2=\pm\Big(\,\begin{matrix}w&-v^{-1}\cr v&0\end{matrix}\,\Big)
\;\;{\rm and}\;\;
c_2=(a_2b_2)^{-1}=\pm\Big(\,\begin{matrix}v^{-2}&0\cr w(v^{-1}-v)&v^2\end{matrix}\,\Big),
\]
where we choose the trace $\pm w$ of $a_2$ and $b_2$ so that they have order $(p+1)/2$. The order of $c_2$, one quarter that of $v$ in ${\mathbb F}_p^*$ and thus dividing $(p-1)/4$, can be chosen to be coprime to $3$ by a suitable choice of $v$, e.g.~taking $v^2=u$ so that $c_2$ has order $2$. By~\cite[Ch.~XII]{Dic} the only maximal subgroups of $G$ containing elements of order $(p+1)/2$ are dihedral groups of order $p+1$, and $a_2$ and $b_2$ cannot both be contained in such a subgroup since they do not commute. It follows that $a_2$ and $b_2$ generate $G$. Since the orders of $a_1, b_1$ and $c_1$ are coprime to those of $a_2, b_2$ and $c_2$, and both triples are hyperbolic, they define a Beauville structure on $G$. The element
\[g=\pm\Big(\,\begin{matrix}0&u\cr u&0\end{matrix}\,\Big)\in G,\]
acting by conjugation, induces an automorphism $\beta_1=\beta_2$ of order $2$ of $G$, transposing the generators $a_i$ and $b_i$ for each $i=1, 2$. Thus $G$ has index $2$ in $N_i=A_i$ for each $i$, with the inclusions of triangle groups as in case~(a). Since $Z(G)=1$ it follows that $N$ is a semidirect product of $D\cong G$ by $\langle\alpha_i\rangle\cong C_2$ (in fact $N\cong G\times C_2$, since $\alpha_i$ induces an inner automorphism $\beta_i$ of $G$), and hence $A^0\cong C_2$. The two triples have inequivalent types, so $A=A^0$.

\subsection{Examples based on Fermat curves}

Let $G=C_t\times C_t$, with two generating triples of type $(t, t, t)$, as in Beauville's original example~\cite{Bea}, where $t=5$, and in Catanese's generalisation~\cite{Cat}, where $t$ is coprime to $6$. Each ${\cal C}_i$ is isomorphic to the Fermat curve ${\cal F}_t$ of genus $(t-1)(t-2)/2$ given by $x^t+y^t+z^t=0$. The triangle group $\Delta_i=\Delta(t, t, t)$ is a normal subgroup of index $6$ in $\tilde\Delta_i=\Delta(2, 3, 2t)$, as in case~(c) in Section~3.1. Since $K_i$ is the commutator subgroup $\Delta_i'$ of $\Delta_i$, it is a characteristic subgroup of $\Delta_i$ and hence normal in $\tilde\Delta_i$, so $A_i=N_i\cong\tilde\Delta_i/K_i$ is an extension of $G\cong\Delta_i/K_i$ by $\tilde\Delta_i/\Delta_i\cong S_3$ for $i=1, 2$. This extension splits, with the normal subgroup given by multiplying the homogeneous coordinates $x, y$ and $z$ by powers of $e^{2\pi i/t}$, and a complement given by permuting them. Since $G$ is abelian, we have $G\times G\leq {\rm ker}\,\theta \leq N$. Thus $M=G\times G$, and $A^0$ contains a normal subgroup $I=M/D=(G\times G)/D\cong G$.

Whether $A^0$ properly contains $I$ depends on the choice of generating triples $a_i, b_i, c_i$ defining the Beauville structure on $G$: specifically, as shown in Section~3.1, we need to decide whether $G$ has automorphisms inducing transpositions or $3$-cycles on both of them. For simplicity, let us take $t$ to be a prime $p\geq 5$, so that two generating triples define a Beauville structure if and only if their images in the projective line ${\mathbb P}^1(p)$, formed by the $1$-dimensional subgroups of $G$, are disjoint. Given a generating triple $a_1, b_1, c_1$ for $G$, the $3$-cycle $(a_1, b_1, c_1)$ extends to a unique automorphism $\beta_1$ of $G$, which decomposes $G\setminus\{1\}$ into $(p^2-1)/3$ cycles $(a_2, b_2, c_2)$ of length $3$. Now $p-1$ of these are scalar multiples of $(a_1, b_1, c_1)$, and if $p\equiv 2$ mod~$(3)$ then each of the remaining $3$-cycles $(a_2, b_2, c_2)$ generates $G$, satisfies $a_2b_2c_2=1$, and has a disjoint image from that of $(a_1, b_1, c_1)$ in ${\mathbb P}^1(p)$, so it gives a Beauville structure on $G$; the automorphism $\beta_2$ of $G$ it induces coincides with $\beta_1$, giving an element $\beta$ of order $3$ in $N$. The situation is similar if $p\equiv 1$ mod~$(3)$, except that in order to generate $G$ the triple $a_2, b_2, c_2$ must now avoid the two $1$-dimensional $\langle\beta_1\rangle$-invariant subgroups of $G$. However, if in either case we choose $a_2, b_2, c_2$ not to form a $3$-cycle of $\beta_1$, as is possible provided $p>5$, then there is no element of order $3$ in $N$.

By contrast, transposing elements of generating triples never induces automorphisms of the Beauville surface $\cal S$. Since $G$ is abelian, a transposition of two elements of a triple induces an automorphism of $G$ fixing the third. Such an automorphism has a $1$-dimensional subgroup of fixed points in $G$, so there cannot be an element of order $2$ in $A^0/I$: the two fixed elements (one from each triple) would be multiples of each other, contradicting condition~(3) for a Beauville structure.

These arguments show that if the triples are chosen to be invariant under the same automorphism $\beta=\beta_1=\beta_2$ of order $3$, then $N$ is a semidirect product of $G\times G$ by $\langle\beta\rangle\cong C_3$, so $A^0$ is a semidirect product of $I\cong G=C_p\times C_p$ by $\langle\beta\rangle$, with $\beta$ acting on $G$ as above. Any other choice of triples (possible if $p>5$) gives $N=G\times G$ and $A^0=I\cong G=C_p\times C_p$. 

We have shown that if $\cal S$ is constructed from the Fermat curve ${\cal F}_t$, where $t$ is a prime $p\geq 5$, then $|{\rm Out}^0{\cal S}|=1$ or $3$. One can extend this result to all $t$ coprime to $6$ by using the natural epimorphism $C_{p^e}\to C_p$ to deal with prime powers, and for general $t$ using the direct product decomposition of $C_t$ based on the prime power factorisation of $t$. See~\cite{GJT} for details, including an enumeration of these Beauville surfaces extending asymptotic estimates by Garion and Penegini~\cite{GP}. We will consider indirect automorphisms in Section~6.

\subsection{A useful construction}

The following lemma will be useful for the next example, and also for later constructions.

\begin{lemma}
Let $G$ be a finite group which is a smooth quotient of $\tilde\Delta:=\Delta(2, 3, n)$ and has no subgroups of index $2$. Then $G$ is also a smooth quotient $\Delta/K$ of $\Delta:=\Delta(t, t, t)$, where $t=n/2$ or $n$ as $n$ is even or odd, with $\tilde\Delta/K\cong G\times S_3$. If $t>3$ then the surface group $K$ has normaliser $N(K)=\tilde\Delta$ in $PSL_2({\mathbb R})$.
\end{lemma}

\noindent{\sl Proof.} First let $n=2t$ be even, so there is a normal surface subgroup $L$ of $\tilde\Delta=\Delta(2, 3, 2t)$ with $\tilde\Delta/L\cong G$. Singerman's normal inclusion~(c) shows that $\tilde\Delta$ has a normal subgroup $\Delta:=\Delta(t, t, t)$ with $\tilde\Delta/\Delta\cong S_3$. Now $\Delta L/\Delta$ is a normal subgroup of $\tilde\Delta/\Delta$, corresponding to a normal subgroup of $S_3$, which must be $1, A_3$ or $S_3$. In the first two cases $\tilde\Delta/\Delta L$ has a quotient isomorphic to $S_3/A_3\cong C_2$, and hence so does $\tilde\Delta/L$, against our hypotheses about $G$. It follows that $\Delta L=\tilde\Delta$, so if we define $K=\Delta\cap L$ then $\Delta/K\cong\tilde\Delta/L\cong G$ and $\tilde\Delta/K\cong G\times S_3$. Since $L$ is a surface group, so is its subgroup $K$, so $G$ is a smooth quotient of $\Delta$. If $t>3$ then the normaliser $N(K)$ of $K$ in $PSL_2({\mathbb R})$ contains $\tilde\Delta$; these are Fuchsian groups, and since Singerman~\cite{Sin2} has shown that $\tilde\Delta$ is maximal among Fuchsian groups, we have $N(K)=\tilde\Delta$.

If $n=t$ is odd we can regard $G$ as a quotient of $\tilde\Delta=\Delta(2, 3, 2t)$ via the natural epimorphism $\tilde\Delta\to\Delta(2, 3, n)$. Although $G$ is no longer a smooth quotient, the only torsion elements in the kernel $L$ are elliptic elements of order $2$, conjugate to $w^t$ where $w$ is the canonical generator of $\tilde\Delta$ of order $2t$. However the image of $w$ in $\tilde\Delta/\Delta$ has order $2$, so only even powers of $w$ lie in $\Delta$; thus $K=\Delta\cap L$ is torsion-free and is therefore a surface group. The rest of the proof is as before.  \hfill$\square$

\begin{cor}
If a finite group $G$ satisfies the conditions of Lemma~4.1 for two mutually coprime values $t_1$ and $t_2$ of $t$, with each $t_i>3$, then $G$ admits a Beauville structure of type $(t_1, t_1, t_1; t_2, t_2, t_2)$ with $A=A^0=I\times S_3\times S_3\cong Z(G)\times S_3\times S_3$.
\end{cor}

\noindent{\sl Proof.} The existence of the Beauville structure follows immediately from Lemma~4.1. We have $N_i\cong G\times S_3$ for $i=1, 2$, with $G$ acting on itself by inner automorphisms, and $S_3$ centralising $G$. The arguments in Section~3.2 then give $N=M\times S_3\times S_2=D\times Z\times S_3\times S_3$, so $A^0=I\times S_3\times S_3\cong Z(G)\times S_3\times S_3$. The two triples have inequivalent types, so by Corollary~3.4 there are no indirect automorphisms. \hfill$\square$

\medskip

As we shall show, this allows the construction of examples in which $|A^0:I|$ attains its upper bound of $36$.

\subsection{Examples with automorphism group $S_3\times S_3$}

As an application of Corollary~4.2, let $G=L_2(p)$ for a prime $p>11$, and let
\[
a_i=\pm\Big(\,\begin{matrix}0&1\cr -1&0\end{matrix}\,\Big),
\quad
b_i=\pm\Big(\,\begin{matrix}u&v\cr w&1-u\end{matrix}\,\Big)
\quad{\rm and}\quad
c_i=(a_ib_i)^{-1}=\pm\Big(\,\begin{matrix}-v&u-1\cr u&w\end{matrix}\,\Big),
\]
where $u(1-u)-vw=1$. Then $a_i$ and $b_i$ have traces $0$ and $\pm 1$, so they have orders $2$ and $3$. For $i=1$ let us take $u=0$, so $w=-1/v$, and choose $v$ to have order $p-1$ in ${\mathbb F}_p^*$, so that $c_1$, having trace $\pm (v+v^{-1})$, has order $(p-1)/2>5$. For $i=2$ let us take $u=1, v=-1$ and $w=1$, so that $c_2$ has order $p$. For each $i$ it follows from Dickson's classification of the maximal subgroups of $G$ in~\cite[Ch.~XII]{Dic} that the triple $a_i, b_i, c_i$ generates $G$. Thus $G$ is a smooth quotient of $\Delta(2, 3, n)$ with $n=(p-1)/2$ and with $n=p$, so it satisfies the conditions of Lemma~4.1 with $t=t_1=(p-1)/4$ or $(p-1)/2$ as $p\equiv \pm 1$ mod~$(4)$, and also with $t=t_2=p$. Since $t_1$ and $t_2$ are coprime, and $Z(G)=1$, Corollaries~3.4 and 4.2 show that $G$ has a Beauville structure of type $(t_1, t_1, t_1; t_2, t_2, t_2)$ with $A=A^0\cong S_3\times S_3$. 

One can construct many similar examples, with $G=L_2(q)$ for suitable prime powers $q$, by using Macbeath's results~\cite{Mac} on generating triples for these groups. One can also use this method to construct examples where $G$ is an alternating group $A_n$. Conder~\cite{Con} has shown that $A_n$ is a quotient of $\Delta(2, 3, 7)$ (i.e.~a Hurwitz group) for all sufficiently large $n$, and Everitt~\cite{Eve}, as part of a more general result on Fuchsian groups, has extended this to all hyperbolic triangle groups, such as $\Delta(2, 3, t)$ for $t\geq 7$. One can ensure that these quotients are smooth (most easily by taking $t$ to be prime), so the preceding arguments show that  $A_n$ has Beauville structures with $A=A^0\cong S_3\times S_3$ for all sufficiently large $n$. Using a direct construction, Fuertes and Gonz\'alez-Diez~\cite[\S5]{FG2} have given an explicit example of such a Beauville structure for $A_{15}$.

\section{Realising abelian groups}

Proposition~3.1 shows that the inner automorphism group $I={\rm Inn}\,{\cal S}$ of a Beauville surface $\cal S$ is a finite abelian group, isomorphic to the centre of the corresponding finite group $G$. The main aim of this section is to show that every finite abelian group can arise as ${\rm Inn}\,{\cal S}$ for some Beauville surface $\cal S$, even if extra restrictions are imposed on ${\rm Out}^0{\cal S}$. The general strategy is as follows. Every finite abelian group $H$ is a direct product of cyclic groups. Special linear groups $SL_d(q)$ have cyclic centres, of all possible orders, so by taking $G$ to be a direct product of suitable special linear groups we can arrange that $Z(G)$ is isomorphic to $H$. Using results of Lucchini~\cite{Luc} we can ensure that these special linear groups are quotients of certain triangle groups, and hence, using a lemma which we shall shortly prove, so is their product $G$. In this way we can construct $G$ to have a Beauville structure, with $I\cong Z(G)\cong H$, and with ${\rm Out}^0{\cal S}$ satisfying various other conditions.

Lucchini's result~\cite{Luc} is as follows:

\begin{prop} For each integer $t\geq 7$ there exists an integer $d_t$ such that $SL_d(q)$ is a quotient of $\Delta(2, 3, t)$ for all $d\geq d_t$ and all prime powers $q$. \hfill$\square$
\end{prop}
 
The centre of $SL_d(q)$ is a cyclic group of order $\gcd(d, q-1)$, consisting of the matrices $\lambda I_d$ where $\lambda^d=1$ in ${\mathbb F}_q$. In order to apply Proposition~5.1 we need the following lemma:

\begin{lemma}
Given any integer $m\geq 1$, there exist infinitely many integers $d\geq 2$, for each of which there are infinitely many prime powers $q$ with $Z(SL_d(q))\cong C_m$.
\end{lemma}

\noindent{\sl Proof.} Let $d=rm$ where $r$ is any integer coprime to $2m$. Since $m$ and $r$ are coprime, it follows from the Chinese Remainder Theorem and Dirichlet's Theorem on primes in arithmetic progressions  that there are infinitely many primes $q$ satisfying $q\equiv 1$ mod~$(m)$ and $q\equiv -1$ mod~$(r)$. In such cases $m$ divides both $d$ and $q-1$, so it divides their highest common factor $h$. Since $d/m=r$, which is an odd divisor of $q+1$ and hence coprime to $q-1$, it follows that $m=h$. Thus $Z(SL_d(q))\cong C_m$. \hfill$\square$

\medskip

Our aim is to construct a group $G$, with $Z(G)$ isomorphic to an arbitrary finite abelian group $H$, by taking a direct product of various groups $G_j$ of type $SL_d(q)$, one for each cyclic direct factor of $H$. We need to show that if each of the groups $G_j$ has a Beauville structure, then so has $G$. In order to achieve this, let us define two groups to be {\sl mutually orthogonal\/} if only the trivial group is a quotient of both of them. This concept was introduced in~\cite{BJNS}, and the following result generalises Corollary~18(a) of that paper:

\begin{lemma}
Let $K_1,\ldots, K_k$ be normal subgroups of a group $\Gamma$ such that the quotient groups $G_j=\Gamma/K_j$ are mutually orthogonal, and let $K=K_1\cap\cdots\cap K_k$. Then $\Gamma/K\cong G_1\times\cdots\times G_k$.
\end{lemma}

\noindent{\sl Proof.} Using induction on $k$, it is sufficient to consider the case $k=2$. In this case $K_1K_2=\Gamma$, for otherwise $\Gamma/K_1K_2$ is a nontrivial common quotient of $G_1$ and $G_2$. Then $\Gamma/K=K_1/K\times K_2/K\cong \Gamma/K_2\times\Gamma/K_1\cong G_2\times G_1\cong G_1\times G_2$. \hfill$\square$

\begin{cor}
Let $G_1,\ldots G_k$ be mutually orthogonal finite groups. If each $G_j$ admits a Beauville structure, then so does $G:=G_1\times\cdots\times G_k$.
\end{cor}

\noindent{\sl Proof.} Let triples $(a_{ij}, b_{ij}, c_{ij})$ for $i=1, 2$ define a Beauville structure on $G_j$ for each $j=1,\ldots, k$. Define elements $a_i, b_i, c_i$ of $G$ by $a_i=(a_{i1},\ldots, a_{ik})$, and so on. Then $a_ib_ic_i=1$ for each $i$, and Beauville condition~(2) is satisfied since it is satisfied in at least one (in fact every) direct factor $G_j$. If some power $a_1^r$ of $a_1$ is conjugate in $G$ to a power $a_2^s$ of $a_2$, then the same applies to $a_{1j}^r$ and $a_{2j}^s$ in each $G_j$, so $a_{1j}^r=1$ for each $j$ and hence $a_1^r=1$; a similar argument applies to any other pair chosen from $a_1, b_1, c_1$ and $a_2, b_2, c_2$. For each $i$, by mapping the two canonical generators of the free group $\Gamma=F_2$ to $a_{ij}$ and $b_{ij}$, we represent each $G_j$ as a quotient $\Gamma/K_j$ of $\Gamma$. Lemma~5.2 then shows that $G$ is also a quotient of $\Gamma$, generated by $a_i$ and $b_i$. Thus both triples generate $G$, giving a Beauville structure on this group. \hfill$\square$

\medskip

Note that if the Beauville structure on each $G_j$ has type $(l_{1j}, m_{1j}, n_{1j}; l_{2j}, m_{2j}, n_{2j})$, then that on $G$ has type $(l_1, m_1, n_1; l_2, m_2, n_2)$ where $l_i={\rm lcm}\,(l_{i1},\ldots, l_{ik})$, etc.

\begin{lemma} Distinct groups $SL_d(q)$ for $d\geq 2$ are mutually orthogonal, except that $SL_2(4)$ and $SL_2(5)$ have a common quotient $SL_2(4)\cong A_5\cong L_2(5)$, and $SL_3(2)$ and $SL_2(7)$ have a common quotient $SL_3(2)\cong L_2(7)$.
\end{lemma}

\noindent{\sl Proof.} Apart from $SL_2(2)$ and $SL_2(3)$, these groups are all perfect. Their only simple quotients are the groups $L_d(q)$, and these are mutually non-isomorphic apart from the isomorphisms $L_2(4)\cong L_2(5)$ and $L_3(2)\cong L_2(7)$ (see~\cite{ATLAS}, for instance). The solvable groups $SL_2(2)$ and $SL_2(3)$ have only $C_2$ and $C_3$ respectively as simple quotients, so this argument extends to them. \hfill$\square$

\begin{thm}
Each finite abelian group $H$ is isomorphic to the inner automorphism group ${\rm Inn}\,{\cal S}$ of a Beauville surface $\cal S$ with ${\rm Aut}\,{\cal S}={\rm Aut}^0{\cal S}\cong H\times S_3\times S_3$.
\end{thm}

\noindent{\sl Proof.} We will apply Corollary~4.2 to a suitable group $G$ with $Z(G)\cong H$. If $H$ is the identity group we can use the example in Section~4.5, so we may assume that $H\cong C_{m_1}\times\cdots\times C_{m_k}$ for integers $m_j\geq 2$. Let us choose two distinct primes $t_1, t_2\geq 7$. By Proposition~5.1, if $d\geq\max\{d_{t_1}, d_{t_2}\}$ then $SL_d(q)$ is a quotient of $\tilde\Delta_i=\Delta(2, 3, t_i)$ for all $q$ and for each $i=1, 2$.  Since $2, 3$ and $t_i$ are primes it must be a smooth quotient. By Lemmas~5.2 and 5.5 we can therefore choose mutually orthogonal groups $G_1,\ldots.,G_k$ of type $SL_d(q)$ which are smooth quotients of $\tilde\Delta_i$ for $i=1, 2$ and have centres $Z(G_j)\cong C_{m_j}$ for $j=1,\ldots, k$. The group $G=G_1\times\cdots\times G_k$ has centre $Z(G)\cong H$, and by Lemma~5.3 it is also a smooth quotient of each $\tilde\Delta_i$. Having no subgroups of index $2$, $G$ satisfies the hypotheses of Lemma~4.1, so by Corollary~4.2 it has a Beauville structure of type $(t_1, t_1, t_1; t_2, t_2, t_2)$ with $I\cong H$ and $A=A^0\cong H\times S_3\times S_3$. \hfill$\square$

\medskip

\noindent{\bf Remarks. 1.} One can obtain the slightly weaker result that every finite abelian group is isomorphic to ${\rm Inn}\,{\cal S}$ for some Beauville surface $\cal S$ by combining Corollary~5.4 with recent results of Fairbairn, Magaard and Parker~\cite{FMP} and of Garion, Larson and Lubotzky~\cite{GLL} which show that, with just finitely many exceptions, the groups $SL_d(q)$ for $d\geq 2$ all admit Beauville structures. This avoids the use of Lemma~5.2, since there is now no requirement that $d$ should be sufficiently large. 

\smallskip

\noindent{\bf 2.} The construction used to prove Theorem~5.6 can be adapted to produce other Beauville surfaces, still realising $H$ as their inner automorphism group, but with smaller direct outer automorphism groups. For instance, let $G_0$ be a symmetric group $S_{11}$, with the Beauville structure constructed by Fuertes and Gonz\'alez-Diez in~\cite{FG} and described in Section~4.1, having type $(2, 20, 9; 2, 10, 11)$ and only the identity automorphism. If $G_1,\ldots, G_k$ are as in the proof of Theorem~5.6, then since $G_0$ is mutually orthogonal to them, Corollary~5.4 implies that the group $G=G_0\times G_1\times\cdots\times G_k$ has a Beauville structure. Since $Z(G_0)=1$, this structure has $I\cong Z(G)\cong H$ as before. However, the type of this new structure consists of two distinct triples, each with no repetitions, so we have $A=A^0=I$. This proves:

\begin{thm}
Every finite abelian group is isomorphic to the automorphism group ${\rm Aut}\,{\cal S}$ of some Beauville surface $\cal S$. \hfill$\square$
\end{thm}

Proposition~3.2 shows that, for any Beauville surface $\cal S$, the group ${\rm Out}^0{\cal S}$ is isomorphic to a subgroup of $S_3\times S_3$. Theorems~5.6 and 5.7 show that the two extreme cases can arise, where ${\rm Out}^0{\cal S}$ is the whole group or the trivial group, with ${\rm Inn}\,{\cal S}$ isomorphic to an arbitrary finite abelian group $H$. Similar constructions, using suitable normal inclusions of $\Delta_i$ in $\tilde\Delta_i$, show that other intermediate subgroups of $S_3\times S_3$ can also arise, specifically from direct automorphism groups of the form $A^0\cong H\times H_1\times H_2$, where each $\tilde\Delta_i/\Delta_i\cong H_i\leq S_3$.

\medskip

\noindent{\bf Example 5.1} Let us take $\Delta_1=\Delta(7, 7, 7)$ and $\tilde\Delta_1=\Delta(2, 3, 14)$ (as in the proof of Theorem~5.6, with $t_1=7$); by Proposition~5.1, if $d\geq d_8$ then $SL_d(q)$ is a quotient of $\Delta_2:=\Delta(2, 3, 8)$ for all $q$, necessarily smooth since $\Delta(2, 3, 4)$ ($\cong S_4$) has only $SL_2(2)$ ($\cong S_3$) as a quotient of this type. Since $7$ is coprime to $2, 3$ and $8$, the construction used for Theorem~5.6 yields Beauville structures of type $(7, 7, 7; 2, 3, 8)$ with $I\cong H$; the normal inclusion $\Delta_1\triangleleft\tilde\Delta_1$ gives $H_1\cong S_3$ as before, while the maximality of the triangle group $\Delta_2$ (see~\cite{Sin2}) implies that $\tilde\Delta_2=\Delta_2$ and hence $H_2=1$. The corresponding Beauville surfaces $\cal S$ therefore have ${\rm Aut}^0{\cal S}\cong H\times S_3$, so ${\rm Out}^0{\cal S}\cong S_3$.

\medskip

\noindent{\bf Example 5.2} A similar construction, taking $\tilde\Delta_2=\Delta(2, 3, 8)$ and $\Delta_2$ its subgroup $\Delta(3, 3, 4)$ of index $2$, gives Beauville structures of type $(7, 7, 7; 3, 3, 4)$ with ${\rm Aut}^0{\cal S}\cong H\times S_3\times C_2$ and thus ${\rm Out}^0{\cal S}\cong S_3\times C_2$.

%%[Can we have a direct factor $H_i\cong C_3$? Can $A^+/I$ be an arbitrary subgroup of $S_3\times S_3$?]

\medskip

In the proofs of Theorem~5.6 and 5.7, and in the above examples, $I$ is in the centre of $A^0$, but Section~4.3 shows that this is not always the case, at least when $I\cong C_p\times C_p$ for some prime $p$. Here we give further examples, with $I$ isomorphic to other groups.

\medskip

\noindent{\bf Example 5.3} Let $S$ be a non-identity finite group with a generating triple $(x, y, z)$ of type $(l, m, n)$. Let $G=S\times S\times S$, and let $G^*$ be the subgroup of  $G$ generated by
\[a_1=(x, y, z), \; b_1=(y, z, x)\quad{\rm and} \quad c_1=(z ,x, y),\]
so that $a_1b_1c_1=1$ since $xyz=yzx=zxy=1$. Suppose that $l, m$ and $n$ are mutually coprime (for instance, $S$ could be a Hurwitz group, with $(l, m, n)=(2, 3, 7)$). Then $G^*$ contains $a_1^{mn}=(x^{mn},1,1)$ and hence contains $(x,1,1)$ since $x$ is a power of $x^{mn}$. Similar arguments show that $G^*$ contains generators for all three direct factors of $G$. Thus $G^*=G$, and hence $G$ is a quotient of $\Delta:=\Delta(t,t,t)$ where $t=lmn$.

Since $x^{-1}z^{-1}y^{-1}=1$, essentially the same argument shows that the elements \[a_2=(x^{-1}, z^{-1}, y^{-1}), \; b_2=(z^{-1}, y^{-1}, x^{-1}) \quad {\rm and} \quad c_2=(y^{-1}, x^{-1}, z^{-1})\]
of order $t$ generate $G$ and satisfy $a_2b_2c_2=1$. No non-identity power of $y$ can be conjugate in $S$ to a power of $z$ (since they have coprime orders), so by considering their second coordinates we see that the same applies to $a_1$ and $a_2$ in $G$. In fact this applies to any pair of elements chosen from the first and the second of these two triples. Since $t>3$ they therefore form a Beauville structure in $G$.

The automorphism $\theta: (g_1, g_2, g_3)\mapsto(g_2, g_3, g_1)$ of $G$ has order $3$ and permutes $a_i, b_i$ and $c_i$ cyclically for each $i$, so the corresponding extension $\tilde G$ of $G$ by $\langle\theta\rangle$ (the wreath product $S\wr C_3$) is a quotient of $\tilde\Delta=\Delta(3,3,t)$, with the same kernel $K$ as $G$. Thus $|N_i:G|$ divisible by $3$ for $i=1, 2$.

To show that $|N_1:G|=3$ it is sufficient to show that $G$ has no automorphism transposing $a_1$ and $b_1$. Such an automorphism would transpose $a_1^{mn}=(x^{mn},1,1)$ and $b_1^{mn}=(1,1,x^{mn})$, and hence $(x,1,1)$ and $(1,1,x)$; similarly it would transpose $a_1^{ln}$ and $b_1^{ln}$, and hence $(1,y,1)$ and $(y,1,1)$; since $(1,1,x)$ and $(1,y,1)$ commute, so must $(x,1,1)$ and $(y,1,1)$, which is impossible since their product $(z^{-1},1,1)$ has order $n$ coprime to $lm$. A similar argument shows that $|N_2:G|=3$. 

Since $N_1$ and $N_2$ induce the same group of automorphisms of $G$ (both acting as $\tilde G$), it follows that the direct automorphism group $A^0$ of the corresponding Beauville surface is a semidirect product of $I\cong Z(G)=Z(S)\times Z(S)\times Z(S)$ by $\langle\theta\rangle\cong C_3$, with $\theta$ permuting the three direct factors $Z(S)$ in a $3$-cycle, so that $A^0\cong Z(S)\wr C_3$. As in the proof of Theorem~5.6, by applying Proposition~5.1 to $\Delta(2, 3, n)$ with $n\;(\geq 7)$ coprime to $6$ one can choose $S$ so that $Z(S)$ is isomorphic to any given finite abelian group.

\section{Beauville surfaces with indirect automorphisms}

We now return to the situation in Section~3.3, where $\cal S$ has an indirect automorphism $\tau:(p_1, p_2)\mapsto (p_2\phi_2, p_1\phi_1)$. Here $\phi_1$ and $\phi_2$ are isomorphisms ${\cal C}_1\to{\cal C}_2$ and ${\cal C}_2\to{\cal C}_1$, and the representations $\rho_i$ of $G$ on the curves ${\cal C}_i$ satisfy $\zeta\circ\rho_1=\phi_2^{-1}\rho_2\phi_2$ and $\zeta\circ\rho_2=\phi_1^{-1}\rho_1\phi_1$ for some $\zeta\in {\rm Aut}\,G$. Any other indirect automorphism $\tau'$ is obtained by composing $\tau$ with a direct automorphism $(\alpha_1, \alpha_2)$; these automorphisms $\alpha_i$ of the curves ${\cal C}_i$ can be absorbed into new isomorphisms $\phi_i'$, leaving $\zeta$ unchanged, so $\zeta$ is independent of the choice of $\tau\in A\setminus A^0$.

\subsection{Normality of the inner automorphism group}

A direct automorphism of $\cal S$ has the form $\alpha=(\alpha_1, \alpha_2)\in N_1\times N_2$, with $\alpha_1$ and $\alpha_2$ inducing the same automorphism of $G$ by conjugation (see Section~3.1). A simple calculation shows that the action of $\tau$ by conjugation on $A^0$ is given by $\alpha^{\tau}=(\alpha_2^{\phi_2}, \alpha_1^{\phi_1})$. Now $\alpha$ is an inner automorphism of $\cal S$ if and only if each $\alpha_i=g_i\rho_i$ for some $g_i\in G$ with $z:=g_1g_2^{-1}\in Z$, in which case $\alpha^{\tau}=((g_2\rho_2)^{\phi_2}, (g_1\rho_1)^{\phi_1})=((g_2\zeta)\rho_1, (g_1\zeta)\rho_2)$ with $(g_2\zeta)(g_1\zeta)^{-1}=z^{-1}\zeta\in Z$. Thus the indirect automorphisms normalise $I$, acting by conjugation as the automorphism
\begin{equation}
z\mapsto z^{\tau}=z^{-1}\zeta=(z\zeta)^{-1},
\end{equation}
where we identify $I$ with $Z$ by means of the isomorphism $\alpha\mapsto z$. Since $I$ is normal in $A^0$, it follows that $I$ is normal in $A$.

As explained at the end of Section~3.3, there is an action of ${\rm Out}^0{\cal S}=A^0/I$ as a subgroup of $S_3\times S_3$, permuting the set $B^2$ where $B=\{0, 1, \infty\}$, and preserving the relations $\equiv_i$ on $B^2$ of having the same $i$th component. Any indirect automorphism $\tau$ acts on $B^2$ by transposing these two equivalence relations, so it induces an element of the wreath product $S_3\wr S_2$, the largest group of permutations of $B^2$ preserving $\{\equiv_1, \equiv_2\}$. This proves the following analogue of Proposition~3.2:

\begin{prop} If $\cal S$ is a Beauville surface, obtained from a Beauville structure on a group $G$, then the automorphism group ${\rm Aut}\,{\cal S}$ of $\cal S$ has an abelian normal subgroup $I\cong Z(G)$ with quotient group ${\rm Out}\,{\cal S}$ isomorphic to a subgroup of $S_3\wr S_2$. In particular ${\rm Aut}\,{\cal S}$ is solvable, of derived length at most $4$, and it has order dividing $72\,|Z(G)|$. \hfill$\square$
\end{prop}

\subsection{Triangle group inclusions}

If $\cal S$ is a Beauville surface in which ${\cal C}_1$ and ${\cal C}_2$ are isomorphic, then we may regard them as a single curve $\cal C$, uniformised by the same surface group $K$. This must be a normal subgroup of hyperbolic triangle groups $\Delta_1$ and $\Delta_2$, with each $\Delta_i/K\cong G$. Since each $\Delta_i$ normalises $K$ it is contained in the normaliser $N(K)$ of $K$ in $PSL_2({\mathbb R})$, and as shown by Singerman~\cite{Sin2}, this must also be a triangle group $\Delta^*$.

If $\cal S$ is to have indirect automorphisms, then $\Delta_1$ and $\Delta_2$ must be of the same type. Now it has been shown by Girondo and Wolfart in~\cite[Theorem 13]{GW} that, except in one special case, if two hyperbolic triangle groups $\Delta_1$ and $\Delta_2$ of the same type are contained in another triangle group $\Delta^*$, then they are conjugate in $\Delta^*$. In our situation, where $\Delta^*=N(K)$, we may without loss of generality conjugate one of them by an element of $\Delta^*$ and thus take them to be the same triangle group $\Delta$. Then the two actions of $\Delta$ on $\cal C$ as $\Delta_1$ and $\Delta_2$,  having the same kernel $K$, must (up to equivalence) differ by an automorphism of $G$; if they are to form a Beauville structure, this must be an outer automorphism.

The exceptional case arises when $\Delta_1$ and $\Delta_2$ both have type $(n, 2n, 2n)$ for some integer $n\geq 3$, and are distinct subgroups of index $2$ in $\Delta^*=\Delta(2, 2n, 2n)$, namely the normal closures in $\Delta^*$ of its two canonical generators of order $2n$; these are clearly not conjugate in $\Delta^*$, though they are conjugate in a triangle group $\Delta(2, 4, 2n)$ which contains $\Delta^*$ as a subgroup of index~$2$. (These inclusions are all of Singerman's type (a) for various choices of the parameters $s$ and $t$.) If this situation is to yield a Beauville surface with indirect automorphisms, then not only must $\Delta_1/K$ and $\Delta_2/K$ be isomorphic (to $G$), but the corresponding two actions of $G$ on $\cal C$ must differ only by automorphisms of $\cal C$ and $G$, as shown in Section~6.1. However, this implies that $\Delta_1$ and $\Delta_2$ must be conjugate in $N(K)=\Delta^*$, which is false, so this case cannot lead to Beauville surfaces with indirect automorphisms. We will therefore assume from now on that $\Delta_1=\Delta_2=\Delta$.

\subsection{Examples of indirect automorphisms}

\medskip

\noindent{\bf Example 6.1} The examples of Beauville surfaces $\cal S$ in~\cite{Bea, Cat}, based on the Fermat curves ${\cal F}_t$, have $\Delta=\Delta(t, t, t)$, $K=\Delta'$, $\Delta^*=\Delta(2, 3, 2t)$ and $G=C_t\times C_t$, with $t$ coprime to $6$; here the inclusion of $\Delta$ in $\Delta^*$ is of Singerman's type~(c). In Section~4.3 we showed that if $t$ is a prime $p\geq 5$ then ${\rm Out}^0{\cal S}$ has order $1$ or $3$; it follows that ${\rm Out}\,{\cal S}$ has order $1, 2, 3$ or $6$ (see~\cite{GJT} for the extension to all $t$ coprime to $6$). Here we give examples of all four cases.

Two generating triples $(a_i, b_i, c_i)$ in $G=C_p\times C_p$ yield a Beauville surface $\cal S$ if and only if their images in the projective line ${\mathbb P}^1(p)$ form disjoint sets $\Sigma_i$. Each $A_i$ is a semidirect product of the abelian group $G_i=G\rho_i\cong G$ by $S_3$, so equivalence of representations $\rho_i$ corresponds to permuting the elements in a generating triple. By Proposition~3.3, if $\cal S$ has an indirect automorphism then some element of $PGL_2(p)$ transposes $\Sigma_1$ and $\Sigma_2$.

We can regard $G$ as a $2$-dimensional vector space over ${\mathbb F}_p$, and without loss of generality we can choose coordinates so that the first generating triple is $a_1=(1,0)$, $b_1=(0,1)$, $c_1=(-1,-1)$. The matrix
\[B=\Big(\,\begin{matrix}0&-1\cr 1&-1\end{matrix}\,\Bigr)\in GL_2(p),\]
acting on the left of column vectors, defines an automorphism $\beta$ of $G$ inducing a $3$-cycle $(a_1, b_1, c_1)$ on this triple. If we let each $(i,j)\neq(0,0)$ correspond to the point $i/j\in{\mathbb P}^1(p)$, then this generating triple corresponds to the triple $\Sigma_1=\{\infty, 0, 1\}$ in ${\mathbb P}^1(p)$. If we choose the second generating triple $a_2, b_2, c_2$ to form a $3$-cycle of $\beta$, inducing a triple $\Sigma_2\subset{\mathbb P}^1(p)$ disjoint from $\Sigma_1$, we obtain a Beauville surface $\cal S$ with ${\rm Out}^0{\cal S}\cong C_3$ (see Section~4.3). Then $\cal S$ also has indirect automorphisms if and only if there is an automorphism $\zeta$ of $G$ transposing these two generating triples. In this case, ${\rm Out}\,{\cal S}$ has order $6$ and is therefore isomorphic to $C_6$ or $S_3$ as $\zeta$ centralises or inverts $\beta$.

Apart from $-I$, which does not lead to a Beauville structure, the involutions in $GL_2(p)$ which commute with $B$ are the matrices of the form
\[Z=\Big(\,\begin{matrix}u&-2u\cr 2u&-u\end{matrix}\,\Bigr)\]
where $3u^2=-1$. By quadratic reciprocity, ${\mathbb F}_p$ contains such elements $u$ if and only if $p\equiv 1$ mod~$(3)$. The automorphism $\zeta$ induced by $Z$ transposes $a_1, b_1, c_1$ with the generating triple $a_2=(u,2u)$, $b_2=(-2u,-u)$, $c_2=(u,-u)$. This corresponds to a triple $\Sigma_2=\{1/2, 2, -1\}\subset{\mathbb P}^1(p)$ disjoint from $\Sigma_1$, so these generating triples $a_i, b_i, c_i$ form a Beauville structure in $G$ with ${\rm Out}\,{\cal S}\cong C_6$.

Similarly, the involutions in $GL_2(p)$ inverting $B$ are the matrices
\[Z=\Big(\,\begin{matrix}u&v\cr u+v&-u\end{matrix}\,\Bigr)\]
where $u^2+uv+v^2=1$. Such a matrix induces an automorphism $\zeta$ of $G$ which permutes the $3$-cycles of $\beta$, reversing their cyclic order. It transposes the generating triple $a_1, b_1, c_1$ with the triple $a_2=(u,u+v)$, $b_2=(v,-u)$, $c_2=(-u-v,-v)$, these triples forming a Beauville structure if and only if $u, v, u+v\neq 0$. We then find that ${\rm Out}\,{\cal S}\cong \langle B, Z\rangle \cong S_3$. Suitable values of $u$ and $v$ exist for all primes $p>7$: in the projective plane ${\mathbb P}^2(p)$, the conic $u^2+uv+v^2=w^2$ has $p+1$ points, of which either two or none are on the line at infinity $w=0$ as $p\equiv 1$ or $-1$ mod~$(p)$, so provided $p>7$ there are points on the affine curve $u^2+uv+v^2=1$ besides the six points $(\pm 1,0), (0,\pm 1)$ and $(\pm 1, \mp 1)$ we need to avoid.

If we form a Beauville structure by choosing for the second triple a $3$-cycle of $\beta$ which is not of one of the above two types, then we obtain a Beauville surface $\cal S$ with ${\rm Out}\,{\cal S}={\rm Out}^0{\cal S}\cong C_3$.

\medskip

\noindent{\bf Example 6.2} Further examples are given by the groups
\[G=\langle g,h\mid g^{p^e}=h^{p^e}=1,\,h^g=h^q\rangle,\]
where $p$ is an odd prime and $q=1+p^f$ with $f=1, 2, \ldots, e$. These groups arose in~\cite{JNS} in connection with regular embeddings of complete bipartite graphs, and various associated regular dessins were studied in~\cite{JSW}. Such a group $G$ is a semidirect product of a normal subgroup $\langle h\rangle$ by $\langle g\rangle$, both cyclic of order $n=p^e$, so each element of $G$ has the unique form $g^ih^j$ where $i, j\in{\mathbb Z}_n$. (This group is a direct product if and only if $e=f$, in which case we have a Fermat curve, as in Section~4.3 and Example~6.1.) As shown in~\cite[Corollary~10]{JNS}, $G$ has exponent $n$, the elements of this order being those of $G\setminus\Phi$, where $\Phi$ is the Frattini subgroup consisting of those elements with $i\equiv j\equiv 0$ mod~$(p)$. By~\cite[Corollary~11]{JNS}, two cyclic subgroups of order $n$ in $G$ have trivial intersection if and only if their images in $G/\Phi\cong C_p\times C_p$ also do. It follows from this that two triples form a Beauville structure in $G$ if and only if their images form a Beauville structure in $G/\Phi$. We therefore obtain Beauville structures in $G$ provided $p\geq 5$. Let us choose $a_1=gh$, $b_1=g^{-4}h^{-2}$ and $c_1=(a_1b_1)^{-1}$, so the corresponding values of $i/j\in{\mathbb P}^1(p)$ are $1, 2$ and $3$. The above presentation shows that $G$ has an automorphism $\zeta$ of order $2$ fixing $g$ and inverting $h$. Applying $\zeta$ to the triple $a_1, b_1, c_1$, we obtain a second triple $a_2, b_2, c_2$ with $i/j=-1, -2$ and $-3$, values which are disjoint from those for the first triple provided $p\geq 7$. In such cases, Proposition~3.3 shows that the corresponding Beauville surface has an indirect automorphism.

\medskip

It is shown in~\cite[\S 4]{JNS} that the group $G$ in Example~6.2 has centre
\[Z(G)=\langle g^{p^{e-f}}, h^{p^{e-f}}\rangle\cong C_{p^f}\times C_{p^f},\]
so the corresponding surface has non-identity inner automorphisms. For later use we will now give a class of examples with an indirect automorphism, but only the identity inner automorphism.

\medskip

\noindent{\bf Example 6.3} Let $P=\langle g, h\mid g^p=h^p=1, gh=hg\rangle\cong C_p\times C_p$, for a prime $p\geq 13$. The triples
\[\overline a_1=g^2h,\; \overline b_1=g^{-6}h^{-2}, \; \overline c_1=g^4h\quad{\rm and}\quad
\overline a_2= gh^2, \; \overline b_2=g^{-2}h^{-6}, \; \overline c_2=gh^4\]
form a Beauville structure, and are transposed by the automorphism $g\mapsto h, h\mapsto g$ of $P$. 

Now let $G=H_1\times H_2\cong H^2$ where
\[H_1=\langle u, g\mid u^q=g^p=1, u^g=u^{\lambda}\rangle\]
and
\[H_2=\langle v, h\mid v^q=h^p=1, v^h=v^{\lambda}\rangle\]
are isomorphic copies of the non-abelian group $H$ of order $pq$ for a prime $q\equiv 1$ mod~$(p)$, with $\lambda^p=1\neq \lambda$ in ${\mathbb Z}_q$. Thus $G$ is a semidirect product of a normal subgroup $Q=\langle u, v\rangle\cong C_q\times C_q$ by $P=\langle g, h\rangle\cong C_p\times C_p$, with $Z(G)=Z(H_1)\times Z(H_2)=1$. We lift the Beauville structure in $P$ to $G$ by defining triples
\[a_1=u^iv^rg^2h,\; b_1=u^jv^sg^{-6}h^{-2}, \; c_1=u^kv^tg^4h\]
and
\[a_2=u^rv^igh^2, \; b_2=u^sv^jg^{-2}h^{-6}, \; c_2=u^tv^kgh^4,\]
transposed by the automorphism $\zeta$ of $G$ which transposes $g$ and $h$, and $u$ and $v$. The condition that each triple should have product $1$ can be written as
\[1=u^ig^2.u^jg^{-6}.u^kg^4=u^i.g^2u^jg^{-2}.g^{-4}u^kg^4=u^{i+j\lambda^{-2}+k\lambda^4}\]
in $H_1$ and
\[1=v^rh.v^sh^{-2}.v^th=v^r.hv^sh^{-1}.h^{-1}v^th=v^{r+s\lambda^{-1}+t\lambda},\]
in $H_2$, that is,
\[{i+j\lambda^{-2}+k\lambda^4=r+s\lambda^{-1}+t\lambda}=0\]
in ${\mathbb Z}_q$.

All six elements in these triples have order $p$, since their images in $P$ act without fixed points on $Q\setminus\{1\}$. It follows that no non-identity power of $a_1, b_1$ or $c_1$ can be conjugate in $G$ to a power of $a_2, b_2$ or $c_2$, since this property holds for their images in $P$.

By its construction, the triple $a_1, b_1, c_1$ maps onto a generating triple for $P$, so it generates a subgroup $G_0$ of $G$ of order divisible by $p^2$. For simplicity, let us take $i=1, j=-\lambda^2, k=0$ and $r=1, s=-\lambda, t=0$, satisfying the above conditions in ${\mathbb Z}_q$, so the triples are
\[a_1=uvg^2h,\; b_1=u^{-\lambda^2}v^{-\lambda}g^{-6}h^{-2}, \; c_1=g^4h\]
and
\[a_2=uvgh^2, \; b_2=u^{-\lambda}v^{-\lambda^2}g^{-2}h^{-6}, \; c_2=gh^4.\]
The projection in $H_1$ of the commutator $[a_1, c_1]$  is
\[[ug^2, g^4]=g^{-2}u^{-1}.g^{-4}.ug^2.g^4=g^{-2}u^{-1}g^2.g^{-6}ug^6
=u^{-\lambda^2+\lambda^6},\]
and its projection in $H_2$ is
\[ [vh, h]=h^{-1}v^{-1}.h^{-1}.vh.h=h^{-1}v^{-1}h.h^{-2}vh^2=v^{-\lambda+\lambda^2},\]
so $G_0$ contains the element
\[[a_1, c_1]=u^{-\lambda^2+\lambda^6}v^{-\lambda+\lambda^2}\]
of $Q$. Conjugating this by $a_1$, we see that $G_0$ also contains the element
\[u^{-\lambda^4+\lambda^8}v^{-\lambda^2+\lambda^3}\]
of $Q$. These two elements generate $Q$ since the determinant
\[(-\lambda^2+\lambda^6)(-\lambda^2+\lambda^3)
-(-\lambda+\lambda^2)(-\lambda^4+\lambda^8)
=-\lambda^4(\lambda-1)^2(1-\lambda^4)\]
is non-zero in ${\mathbb Z}_q$. Thus $G_0=G$, so the triple $a_1, b_1, c_1$ generates $G$, and hence so does $a_2, b_2, c_2$. These triples therefore form a Beauville structure in $G$. Since $Z(G)=1$ the corresponding Beauville surface $\cal S$ has only the identity inner automorphism, that is, $I=1$. As shown in Example~4.3, the Beauville surface $\overline{\cal S}$ corresponding to $P$ cannot have a direct outer automorphism of order $2$, and it is easy to check that an automorphism of $P$ inducing a $3$-cycle on $\overline a_1, \overline b_1, \overline c_1$ does not leave the triple $\overline a_2, \overline b_2, \overline c_2$ invariant; thus $P$ has no direct outer automorphisms, and hence the same applies to $\cal S$, giving $A^0=1$. However, $\cal S$ has an indirect automorphism of order $2$ induced by the automorphism $\zeta$ of $G$ transposing the two triples $a_i, b_i, c_i$, so $A\cong C_2$.

\subsection{Realising abelian groups}

Using Example~6.3, we will show that every finite abelian group can arise as the inner automorphism group of a Beauville surface admitting indirect automorphisms. First we need the following lemma:

\begin{lemma}
Let $G_0, G_1, \ldots, G_k$ be mutually orthogonal finite groups. For each $j=0, \ldots, k$ let $G_j$ have generating triples $(a_{1j}, b_{1j}, c_{1j})$ and $(a_{2j}, b_{2j}, c_{2j})$ such that
\begin{itemize}
\item $(a_{10}, b_{10}, c_{10})$ and $(a_{20}, b_{20}, c_{20})$ form a Beauville structure for $G_0$, and
\item each $(a_{1j}, b_{1j} c_{1j})$ has type $(l_j, m_j, n_j)$, where $l_j$ divides $l_0$, $m_j$ divides $m_0$ and $n_j$ divides $n_0$ for each $j$.
\end{itemize}
Then the elements $a_1=(a_{10},\ldots, a_{1k}), \ldots, c_2=(c_{20},\ldots, c_{2k})$ form a Beauville structure for the group $G=G_0\times\cdots\times G_k$.
\end{lemma}

(Note that for $j=1, \ldots, k$ we do not require the two triples in $G_j$ to form a Beauville  structure. Indeed, in some cases we can (and will) take $(a_{1j}, b_{1j}, c_{1j})$ and  $(a_{2j}, b_{2j}, c_{2j})$ to be the same triple.)

\medskip

\noindent{\sl Proof.}  If some power $a_1^r$ of $a_1$ is conjugate in $G$ to a power $a_2^s$ of $a_2$, then the same applies to their projections $a_{10}^r$ and $a_{20}^s$ in $G_0$, so $a_{10}^r=1$ since the two triples in $G_0$ form a Beauville structure. Thus $r$ is divisible by $l_0$, and hence by $l_j$ for every $j$, so $a_{1j}^r=1$ and hence $a_1^r=1$. A similar argument applies to any other pair chosen from the triple $(a_1, b_1, c_1)$ and the triple $(a_2, b_2, c_2)$. The rest of the proof follows that for Lemma~5.3. \hfill$\square$

\medskip

Given an abelian group $H$, the {\sl generalised dihedral group\/} ${\rm Dih}\,H$ is the semidirect product of $H$ by a complement $C_2$ inverting $H$ by conjugation. The next result is an analogue of Theorem~5.7:

\begin{thm}
Each finite generalised dihedral group is isomorphic to ${\rm Aut}\,{\cal S}$ for some Beauville surface $\cal S$ with an indirect automorphism.
\end{thm}

\noindent{\sl Proof.} Let $G_0$ be the group of order $p^2q^2$ denoted by $G$ in Example~6.3, with a Beauville structure in which two triples $(a_{10}, b_{10}, c_{10})$ and $(a_{20}, b_{20}, c_{20})$ of type $(p, p, p)$ are transposed by an automorphism of $G_0$, for some prime $p\geq 13$ (for instance, we could take $p=13$ and $q=53$). Given any finite abelian group $H\cong C_{m_1}\times\cdots\times C_{m_k}$, we can use Proposition~5.1 and Lemma~4.1 to choose mutually orthogonal groups $G_j\;(j=1,\ldots, k)$ of type $SL_d(q)$, so that each $G_j$ is a smooth quotient of $\Delta(2, 3, p)$, and hence of $\Delta(p, p, p)$, with $Z(G_j)\cong C_{m_j}$. In Lemma~6.2 we take $(a_{1j}, b_{1j}, c_{1j})=(a_{2j}, b_{2j}, c_{2j})$ to be the corresponding generating triple of $G_j$ of type $(p, p, p)$ for $j=1,\ldots, k$. By Lemma~6.2 the elements $a_1=(a_{10},\ldots, a_{1k}), \ldots, c_2=(c_{20},\ldots, c_{2k})$ form a Beauville structure for $G:=G_0\times G_1\times\cdots\times G_k$. Since $Z(G_0)=1$ the corresponding Beauville surface $\cal S$ has inner automorphism group $I\cong Z(G)\cong H$. The two triples $(a_1, b_1, c_1)$ and $(a_2, b_2, c_2)$ are transposed by an outer automorphism of $G$ which acts as in Example~6.3 on $G_0$, and as the identity on $G_j$ for each $j=1,\ldots, k$, thus inducing an indirect automorphism $\tau$ of $\cal S$.

It is shown in Example~6.3 that the Beauville surface constructed there from $G_0$ has no direct outer automorphisms, so the same applies to $\cal S$. Thus ${\rm Aut}\,{\cal S}$ is a semidirect product of $I$ and a group $\langle\tau\rangle\cong C_2$. Now $\tau$ acts by conjugation on $I$ as $(\alpha_1, \alpha_2)\mapsto(\alpha_2,\alpha_1)$, and since we identify $(\alpha_1,\alpha_2)$ with the element $\alpha_1\alpha_2^{-1}\in Z(G)$ it follows that $\tau$ acts on $I$ by inverting each element, so ${\rm Aut}\,{\cal S}\cong {\rm Dih}\,I\cong {\rm Dih}\,H$.
\hfill$\square$ 

\medskip

In Theorem~6.3, the order of the outer automorphism group ${\rm Out}\,{\cal S}$ attains its lower bound among all Beauville surfaces with indirect automorphisms, namely $|{\rm Out}\,{\cal S}|=2$. We will now construct examples in which it attains its upper bound of $72$.

\medskip

\noindent{\bf Example 6.4} Let $G=G_1\times G_2$, where $G_1=L_2(5^2)$ and $G_2=L_2(3^3)$. Our aim is to construct a Beauville structure of type $(13, 13, 13; 13, 13, 13)$ for $G$, even though it follows easily from Sylow's Theorems that neither $G_1$ nor $G_2$ can have a structure of this type, since each has cyclic Sylow $13$-subgroups.

It can be seen from their character tables and lists of maximal subgroups in~\cite{ATLAS}, or alternatively deduced from results of Macbeath~\cite{Mac}, that each $G_j$ is a smooth quotient of $\Delta(2, 3, 13)$, and hence by Lemma~4.1 of $\Delta(13, 13, 13)$, giving a generating triple $(a_{1j}, b_{1j}, c_{1j})$ of type $(13, 13, 13)$ with all three generators in the same conjugacy class. In each case, this can be chosen to be any of the six conjugacy classes of elements of order $13$ in $G_j$. Note that the elements of each such class are conjugate to their inverses, but to no other proper powers of themselves.

In $G_1$ we take $a_{21}, b_{21}$ and $c_{21}$ to be the images of $a_{11}, b_{11}$ and $c_{11}$ under the automorphism of order $2$ induced by the Frobenius automorphism $z\mapsto z^5$ of the underlying field ${\mathbb F}_{5^2}$. These elements lie in a conjugacy class consisting of the $5$th powers of those in the class containing $a_{11}, b_{11}$ and $c_{11}$, so an element of the first triple is conjugate to the $r$th power of an element of the second triple if and only if $r\equiv\pm 5$ mod~$(13)$. In $G_2$ we take $a_{22}=a_{12}$, $b_{22}=b_{12}$ and $c_{22}=c_{12}$, so in this case an element of the first triple is conjugate to the $r$th power of an element of the second triple if and only if $r\equiv\pm 1$ mod~$(13)$. It follows that in $G$, no element of the triple $a_1=(a_{11}, a_{12}), b_1=(b_{11}, b_{12}), c_1=(c_{11}, c_{12})$ can be conjugate to a power of an element of the triple $a_2=(a_{21}, a_{22}), b_2=(b_{21}, b_{22}), c_2=(c_{21}, c_{22})$.

Since $G_1$ and $G_2$ are mutually orthogonal (as non-isomorphic simple groups), each of these triples generates $G$, so they form a Beauville structure of type $(13, 13, 13; 13, 13, 13)$ in $G$. The two triples are transposed by an outer automorphism $\zeta$ of $G$, which acts as the field automorphism on $G_1$ and the identity on $G_2$, so by Proposition~3.3 the corresponding Beauville surface $\cal S$ has an indirect automorphism $\tau$. Now ${\rm Inn}\,{\cal S}\cong Z(G)=1$, and Lemma~4.1 implies that each curve ${\cal C}_i\;(i=1, 2)$ has an automorphism group $S_3$ commuting with $G$, so $\cal S$ has direct automorphism group ${\rm Aut}^0{\cal S}\cong S_3\times S_3$. The existence of $\tau$ shows that ${\rm Aut}\,{\cal S}$ properly contains ${\rm Aut}^0{\cal S}$, so by Proposition~6.1 it is isomorphic to $S_3\wr S_2$.

\medskip

We can extend the above example by proving an analogue of Theorem~6.3, in which $|{\rm Out}\,{\cal S}|$ now attains its upper bound of $72$. We can again use Lemma~6.2, but instead of taking $G_0$ to be the group in Example~6.3, as used in the proof of Theorem~6.3, we use the group $L_2(5^2)\times L_2(3^3)$ in Example~6.4. Otherwise, the proof follows that of Theorem~6.3, with $p=13$, except that now {\sl every\/} group $G_j\;(j=0,\ldots, k)$ arises as a quotient of $ \Delta(2, 3, 13)$; this implies that the curves ${\cal C}_1$ and ${\cal C}_2$ each have an automorphism group $S_3$ commuting with $G$, so that ${\rm Out}^0{\cal S}\cong S_3\times S_3$. An indirect automorphism $\tau$, acting as the field automorphism of $L_2(5^2)$, and as the identity on $L_2(3^3)$ and the special linear groups $G_1,\ldots, G_k$, transposes the two direct factors $S_3$ of ${\rm Out}^0{\cal S}$, so ${\rm Out}\,{\cal S}\cong S_3\wr S_2$. This proves:

\begin{thm}
Every finite abelian group is isomorphic to the inner automorphism group of some Beauville surface with outer automorphism group isomorphic to $S_3\wr S_2$. \hfill$\square$
\end{thm}

\end{document}